\documentclass[11pt]{amsart}
\usepackage{amssymb}
\usepackage{verbatim}
\usepackage{amscd}
\usepackage{amsmath}

\thanks{}
\subjclass{} \keywords{}

\newtheorem{theorem}{Theorem}[section]
\newtheorem{corollary}[theorem]{Corollary}
\newtheorem{lemma}[theorem]{Lemma}
\newtheorem{proposition}[theorem]{Proposition}

\theoremstyle{definition}

\numberwithin{equation}{section}

\def\boxit#1{\vbox{\hrule\hbox{\vrule\kern3pt
     \vbox{\kern3pt#1\kern3pt}\kern3pt\vrule}\hrule}}

\newfont{\msam}{msam10}            
\newfont{\msym}{msbm10 scaled\magstep1}            

\newfont{\gotic}{eufm10 scaled\magstep1}

\newcommand{\ra}{\rightarrow}
\newcommand{\lra}{\mbox{\Huge $\longrightarrow$}}

\newcommand{\GL}{{\rm GL}}
\newcommand{\son}{{S^{2n,1}}}
\newcommand{\kra}{\kern-7pt\rightarrow\kern-7pt}

\def\bfv{\mathbf{v}}
\def\bfw{\mathbf{w}}

\def\bfr{\mathbf{r}}

\def\bbr{{\mathbb R}}

\def\bbc{{\mathbb C}}

\def\ra{\rightarrow}

\def\x{\times}

\def\im{\mathrm{Im}}

\def\aff{\mathrm{Aff}}

\def\gl{\mathrm{GL}}

\def\inv{^{-1}}

\def\bz{{\mathbf z}}
\def\lra{\longrightarrow}

\def\ol{\overline}
\def\calz{\mathcal Z}

\def\fm{\phantom{-}}
\numberwithin{equation}{section}

\newcommand{\angles}[1]{{\langle #1 \rangle}}
\def\begtab{\begin{tabbing} WW\=23/02: \= point 1\kill}

\def\NIL1{{\mathcal H^{3}}}
\def\NILn{{\mathcal H^{2n+1}}}

\def\wt{\widetilde}

\def\n2{\mathfrak{N}_2}

\def\frakm{{\mathfrak m}}
\def\fraku{{\mathfrak u}}
\def\frakh{{\mathfrak h}}
\def\frakg{{\mathfrak g}}
\def\frakZ{{\mathfrak Z}}

\begin{document}
\title[HOLONOMY DISPLACEMENTS IN  HOPF BUNDLES]
{Holonomy Displacements in Hopf Bundles over complex hyperbolic
space and the complex Heisenberg groups}
\author{Younggi Choi and Kyung Bai Lee}
\address{ Department of Mathematics Education, Seoul National University,
Seoul 151-748, Korea } \email{yochoi@snu.ac.kr}
\address{Department of Mathematics, University of Oklahoma, Norman, OK 73019, U.S.A.}
\email{kb\_lee@math.ou.edu} \subjclass[2000]{53C29, 53C30, 32Q45}
\keywords{Holonomy displacement; complex hyperbolic space; complex
Heisenberg group} \abstract
 For the ``Hopf bundle'' $S^1\ra \son\ra \bbc H^n$, horizontal lifts of simple closed
curves are studied. Let $\gamma$ be a piecewise smooth, simple
closed curve on a complete
 totally geodesic surface $S$ in the base space. Then the holonomy
displacement along $\gamma$ is given by
$$
V(\gamma)=e^{\lambda \, A(\gamma) i}
$$
where $A(\gamma)$ is the area of the region on the surface $S$
surrounded by $\gamma$; $\lambda=1/2 $ or 0 depending on whether
$S$ is a complex submanifold or not.

We also carry out a similar investigation for the complex
Heisenberg group $\bbr\ra\NILn\ra\bbc^n$.
\endabstract
\maketitle

\section{Introduction}
Consider the Hopf fibration $S^1\ra S^3\ra S^2$. Let $\gamma$ be a
simple closed curve on $S^2$. Pick a point in $S^3$ over
$\gamma(0)$, and take the unique horizontal lift $\wt\gamma$ of
$\gamma$. Since $\gamma(1)=\gamma(0)$, $\wt\gamma(1)$ lies in the
same fiber as $\wt\gamma(0)$ does. We are interested in
understanding the difference between $\wt\gamma(0)$ and
$\wt\gamma(1)$. The following equality was already known
\cite{Pin}:
 $$
V(\gamma)=e^{\frac{1}{2}  A(\gamma) i}
$$
where $V(\gamma)$ is  the holonomy displacement along $\gamma$,
and $A(\gamma)$ is the area of the region surrounded by $\gamma$.

In this paper, we shall investigate this fact to (higher
dimensional) pseudo-spheres and the complex Heisenberg group.
First we look at the fibration of the pseudo-sphere $S^{2n,1}$
$$
\begin{array}{lll}
&S^1\ra S^{2n,1}  \ra \bbc H^n
\end{array}
$$
a principal $S^1$-bundle over the complex hyperbolic space $\bbc
H^n$. Let $S$ be a complete totally geodesic surface in the base
space $\bbc H^n$, and $\xi_S$ be the pullback bundle over $S$. Let
$\gamma$ be a piecewise smooth, simple closed curve on $S$
parametrized by $0\leq t\leq 1$, and $\wt\gamma$ its horizontal
lift. The pullback over the curve $\gamma$ is called a \emph{Hopf
torus} (so $\wt\gamma$ is a curve on the Hopf torus). Then
$$
\wt\gamma(1)= e^{\frac{1}{2} A(\gamma) i}\cdot\wt\gamma(0)
\text{\hskip12pt or\hskip12pt } \wt\gamma(0),
\smallskip
$$
depending on whether $S$ is a complex submanifold or not, where
$A(\gamma)$ is the area of the region on the surface $S$
surrounded by $\gamma$. See {\rm Theorem  \ref{thm-hyper}}.

We also carry out a similar investigation for the complex
Heisenberg group.
 Let $1\ra\bbr\ra\NILn\ra\bbc^n\ra 1$
be the central short exact sequence of the complex Heisenberg
group. Let $S$ be a complete totally geodesic plane in $\bbc^n$,
and $\xi_S$ be the pullback bundle over $S$. Let $\gamma$ be a
piecewise smooth, simple closed curve on $S$. Then
$$
V(\gamma)= e(\xi_S)\cdot A(\gamma)
$$
where $A(\gamma)$ is the area of the region on the surface $S$
surrounded by $\gamma$, and the number $e(\xi_S)$ is determined by
the equality $[\bfv,\bfw]=e(\xi_S) e_{2n+1}$ for an orthonormal
basis $\{\bfv,\bfw\}$ for the tangent space of $S$. See {\rm
Theorem \ref{thm-nil}}.

\section{Preliminaries}
The proof of the statement in the introduction (in the case of the
Hopf fibration $S^1\ra S^3\ra S^2$) uses the Gauss-Bonnet
theorem. For $S^1\ra S^{2,1}\ra \bbc H^1$, such is not available
because the base space is not compact. Therefore, we cannot apply
the arguments in \cite{Pin} directly, and need to develop a new
method of proof. It turns out that $S^{2,1}$ is the building
blocks for higher dimensional cases.

Let $F\ra E\stackrel{p}\ra B$ be a principal $F$-bundle
$(F=\bbr^1\text{ or } S^1)$ of Riemannian manifolds, with $B$ a
2-dimensional complete manifold and $p$ a Riemannian submersion.
For a simple closed cure $\gamma(t), 0\leq t\leq 1$ on $B$, the
\emph{holonomy displacement} $V(\gamma)$ along $\gamma$ is defined
as follows: Let $\wt\gamma(t)$ be the horizontal lift of $\gamma$.
Then
$$
\wt\gamma(1)=V(\gamma)\cdot \wt\gamma(0)
$$
for some $V(\gamma)\in F$. We shall establish a technical lemma
which will be used later.

\begin{lemma}
\label{rect-to-gen} Suppose $V(\gamma)=\lambda(\gamma)
A(\gamma)\,(F=\bbr^1) , \text{\ or\,\, } e^{\lambda(\gamma)
A(\gamma) i}\,(F=S^1 )$ for a constant $\lambda(\gamma)$, where
$A(\gamma)$ is the area of the region on $B$ surrounded by  a
piecewise smooth simple closed curve $\gamma$.
If $\lambda(\gamma)$ is constant for all $\gamma$'s which are the
boundaries of rectangular regions, then it is constant for every
piecewise smooth simple closed curve $\gamma$.
\end{lemma}

\begin{proof}
Let us assume that $F=\bbr$. The case of  $F=S^1$ will be similar.
Let $\gamma_0$ be a curve on $B$. Since the region surrounded by
$\gamma_0$ is compact, we may assume that this region is contained
completely in one local patch. Let
$$
\varphi:\bbr^2\ra U\subset B
$$
be a local chart, and $p\inv(U)\approx U\x F$. For notational
simplicity, we shall identify $\bbr^2$ with $U$ (and suppress
$\varphi$). Let $\Omega(U,\gamma_0(0))$ and $\Omega(\bbr,0)$ be
the space of paths emanating from $\gamma_0(0)$ and $0\in \bbr$,
respectively. For each $\gamma\in \Omega(U,\gamma_0(0))$, let
$\omega_\gamma$ be the unique curve in $\bbr$ so that $
\eta(t)=\gamma(t)\cdot\omega_\gamma(t) $ is the horizontal lift of
$\gamma$. This defines a map
$$
\frakZ : \Omega(U)\lra \Omega(\bbr)
$$
by $\frakZ(\gamma)(t)=\omega_\gamma(t)$. We use the $\sup$ metrics
$ \rho$  on both $\Omega(U)$ and $\Omega(\bbr)$. That is,
$$
\rho (\gamma_1,\gamma_2)=\sup_{t\in [0,1]}
d(\gamma_1(t),\gamma_2(t))
$$
\noindent where $d$ is the distance function on $U$. A similar
definition holds  for $\Omega(\bbr)$. We wish to show that
$\frakZ$ is continuous at $\gamma_0$. Let $\epsilon>0$ be given.
By the continuity of the connection, for each $t\in [0,1]$, there
is an open neighborhood $W$ of $\gamma(t)$ such that
any piecewise smooth curve in $W$ has a horizontal lift which lies
in $W\x (-\epsilon/2,\epsilon/2)$. Since $\gamma_0(I)$ is compact,
we can find $\delta>0$  such that if $\rho
(\gamma_0,\gamma)<\delta$, then $ \rho
(\omega_{\gamma_0},\omega_{\gamma})<\epsilon$. This proves that
$\frakZ$ is continuous.

Any piecewise smooth simple closed curve can be approximated by a
sequence of piecewise linear curves which are sums of boundaries
of rectangular regions. Since $\lambda(\gamma)$ is constant for
rectangular regions, the same is true for any  piecewise smooth
simple closed  curve.
\end{proof}

Next, we need to know all complete totally geodesic submanifolds
of the base space of the principal bundle $S^1\ra \son\ra \bbc
H^n$. Since $\son$ is a symmetric space, the following gives a
complete answer.

\begin{proposition}\cite[XI Theorem 4.3]{KN}
\label{tot-geod-prop} Let $(G,H, \sigma)$ be a symmetric space and
$\frakg=\frakh + \frakm$ the canonical decomposition. Then there
is a natural one-to-one correspondence between the set of linear
subspaces $\frakm'$ of $\frakm$ such that
$[[\frakm',\frakm'],\frakm']\subset \frakm'$ and the set of
complete totally geodesic submanifolds $M'$ through the origin $0$
of the affine symmetric space $M=G/H$, the correspondence being
given by $\frakm'=T_{0}(M')$.
\end{proposition}

\bigskip

\section{The bundle $S^1\ra \son\ra \bbc H^n$}

We shall study the  bundle
$$
U(1)\lra U(1,n)/U(n) \stackrel{p}\lra U(1,n)/(U(1)\x U(n)).
$$
Note that $ U(1,n)/U(n) \cong \son,$ and $ U(1,n)/(U(1)\x U(n))
\cong \bbc H^{n} $ where $\son=H^{1,2n} = \{ ( z_0, \ldots , z_{n}
)  \in \bbc^{n} : \  - |z_{0}|^{2}  +  \sum_{i=1}^{n} | z_{i}
|^{2} = -1 \} $. For more information for  $\son$, see \cite{W}.
We  first consider the case when $n=1$. Rather than using $U(1)\ra
U(1,1)/U(1)\ra U(1,1)/(U(1)\x U(1))$, we shall use
$$U(1)\ra SU(1,1)\ra SU(1,1)/U(1). $$  Here
$ SU(1,1)=\{A\in\GL(2,\bbc)\ :\ AJA^*=J\text{ and } \det(A)=1\} $
where ${\smaller
J=\left[\begin{array}{rr}-1&0\\0&1\end{array}\right] }$.

From now on, we shall use the convention of
$\frak{gl}(n,\bbc)\subset\frak{gl}(2n,\bbr)$ by
$$
\left[\begin{array}{cccc}
z_{11} & z_{12}\\
z_{21} & z_{22}\\
\end{array}\right]
\lra \left[\begin{array}{rrrrrrrr}
x_{11} &-y_{11} &x_{12} &-y_{12}\\
y_{11} & x_{11} &y_{12} & x_{12}\\
x_{21} &-y_{21} &x_{22} &-y_{22}\\
y_{21} & x_{21} &y_{22} & x_{22}\\
\end{array}\right].$$

The group $SU(1,1)$ has the following natural representation into
$\GL(4,\bbr)$:
$${\smaller
w= \left[
\begin{array}{rrrr}
\fm w_1 &w_2 & w_3 & w_4\\
-w_2 &w_1 &-w_4 &w_3\\
w_3 &-w_4 &w_1 &-w_2\\
w_4 &w_3 &w_2 &w_1
\end{array}
\right]}
$$
with the condition $w_1^2 + w_2^2 - w_3^2 - w_4^2 =1$. In fact,
the map
$$
w_1 + w_2 i + w_3 j + w_4 k \longmapsto w
$$
is  a monomorphism  from the unit quaternions into $\gl(4,\bbr)$.
Therefore, $SU(1,1)\cong S^{2,1}$. The circle group
\begin{align*}
S^1&= \left\{ \left[
\begin{array}{ll}
e^{iz}  &0  \\
0& e^{-iz}  \\
\end{array}\right]\ :\ 0\leq z\leq 2\pi
\right\}
\end{align*}
is a subgroup of $SU(1,1)$, and acts on $SU(1,1)$ as right
translations, freely with quotient $\bbc H^1$, the complex
hyperbolic line, giving rise to the fibration
$$
S^1\lra SU(1,1)\lra \bbc H^1.
$$

In order to understand the projection map better, let $\wt w$ be
the ``$i$-conjugate'' of $w$ (replace $w_2$ by $-w_2$). That is,
$${\smaller
\wt w= \left[
\begin{array}{rrrr}
\fm w_1 &-w_2 & w_3 & w_4\\
w_2 &w_1 &-w_4 &w_3\\
w_3 &-w_4 &w_1 &w_2\\
w_4 &w_3 &-w_2 &w_1
\end{array}
\right]. }
$$
Then,
$$
w \wt w= {\tiny \left[\begin{array}{rrrrrrrr}
w_1^2+w_2^2+w_3^2+w_4^2 &0 &2 (w_1 w_3- w_2 w_4) & 2 (w_2 w_3+w_1 w_4)\\
0 &w_1^2+w_2^2+w_3^2+w_4^2 &-2 (w_2 w_3+w_1 w_4) & 2( w_1 w_3- w_2 w_4) \\
2 (w_1 w_3- w_2 w_4) &-2 (w_2 w_3+w_1 w_4) &w_1^2+w_2^2+w_3^2+w_4^2 &0\\
2 (w_2 w_3+w_1 w_4) &2( w_1 w_3- w_2 w_4) &0 &w_1^2+w_2^2+w_3^2+w_4^2\\
\end{array}\right]}
$$
and
 $$
(w_1^2+w_2^2+w_3^2+w_4^2)^2-(2 w_1 w_3-2 w_2 w_4)^2 -(2 w_2 w_3+2
w_1 w_4)^2=1.
$$
Clearly,  $\bbc H^1$  can be identified with the following
 $$
\bbc H^1 = \left\{ \left[
\begin{array}{rrrr}
x & 0 & y & z  \\
0 & x &-z & y\\
y &-z & x & 0\\
z & y & 0 & x
\end{array}\right]
:\ x^2 - y^2 - z^2=1, \ x>0 \right\}.
$$
Therefore, the map
$$
p: SU(1,1)\lra \bbc H^1
$$
defined by $ p(w)=w \wt w $ has the following properties:
\begin{align*}
p(wv)&=w p(v) \wt w\quad\text{for all } w,v\in SU(1,1), \\
p(wv)&=p(w) \quad\text{if and only if}\quad v\in S^1.
\end{align*}
This shows that the map $p$ is, indeed, the orbit map of the
principal bundle $S^1\lra SU(1,1)\lra \bbc H^1$. The Lie group
$SU(1,1)$ will have a left-invariant Riemannian metric given by
the following orthonormal basis on the Lie algebra
${\mathfrak{su}}(1,1)$
\[ {\smaller
e_{1}= \left( \begin{array}{rrrr}
0 & 0 & 1& 0  \\
0 & 0 & 0  & 1\\
1 & 0 & 0 &0\\
0 & 1 & 0 &0
\end{array}\right),\;\;
e_{2}= \left( \begin{array}{rrrr}
0 & 0 & 0& 1  \\
0 & 0 & -1  & 0\\
0 & -1 & 0 &0\\
1 & 0 & 0 &0
\end{array}\right),\;\;
e_{3}= \left( \begin{array}{rrrr}
0 & 1 & 0& 0  \\
-1 & 0 & 0  & 0\\
0 & 0 & 0 &-1\\
0 & 0 & 1 &0\\
\end{array}\right). }
\]

Notice that $e_1$ and $e_2$ correspond to $\left[\begin{array}{ll}
0\\ 1\end{array}\right]$ and $\left[\begin{array}{ll} 0\\
i\end{array}\right]$
 in $\frak{gl}(2,\bbc)$
and $[e_1,e_2]=-2 e_3$. Consider the subset of $SU(1,1)$:
{\smaller
\begin{align*}
T&= \left\{ \left[\begin{array}{ll}
\cosh x &(\sinh x)e^{-i y}\\
(\sinh x)e^{i y} &\cosh x\\
\end{array}\right]\ :\  x\geq 0,\ 0\leq y\leq 2\pi
\right\}\\
&= \left\{ \left[\begin{array}{rrrrrrrr}
\cosh x &0          & (\sinh x)(\cos y) & (\sinh x)(\sin y)\\
0          &\cosh x &-(\sinh x)(\sin y) & (\sinh x)(\cos y)\\
 (\sinh x)(\cos y) &-(\sinh x)(\sin y) &\cosh x &0         \\
 (\sinh x)(\sin y) & (\sinh x)(\cos y) &0          &\cosh x\\
\end{array}\right]
\right\}
\end{align*} }
which is the exponential image of
$$
\frakm= \left\{\left[\begin{array}{cc}
 0 &  \bar{\xi}^{t} \\ \xi &0
\end{array} \right]\ :\ \xi  \in  {\bbc} \right\}.
$$
Of course, $SU(1,1)$ is topologically a product $S^1\x\bbc H^1$.
The map $p$ restricted to $T$ is just the squaring map; that is,
$$
p(w)=w^2,\quad w\in T.
$$
\smallskip

\begin{theorem}
\label{area-u11} Let $S^1\ra SU(1,1)\ra \bbc H^1$ be the natural
fibration. Let $\gamma$ be a piecewise smooth, simple closed curve
on $\bbc H^1$. Then the holonomy displacement along $\gamma$ is
given by
$$
V(\gamma)=e^{\frac{1}{2} A(\gamma) i}\ \in S^1
$$
where $A(\gamma)$ is the area of the region on $\bbc H^1$ enclosed
by $\gamma$.
\end{theorem}

\begin{proof}
Let $\gamma(t)$ be a closed loop on $\bbc H^1$ with
$\gamma(0)=p(I_4)$. Therefore, {\tiny
\begin{align*}
\gamma(t)= \left[\begin{array}{rrrrrrrr}
\cosh 2x(t) &0          &\sinh 2x(t)\cos y(t) &\sinh 2x(t)\sin y(t)\\
0          &\cosh 2x(t) &- \sinh 2x(t)\sin y(t) &\sinh 2x(t)\cos y(t)\\
 \sinh 2x(t)\cos y(t) &-\sinh 2x(t)\sin y(t) &\cosh 2x(t) &0         \\
 \sinh 2x(t)\sin y(t) & \sinh 2x(t)\cos y(t) &0          &\cosh 2x(t)\\
\end{array}\right].
\end{align*}}
\par\noindent
Let
\begin{align*}
{\tiny \wt\gamma(t)= \left[\begin{array}{rrrrrrrr}
\cosh x(t) &0          &\sinh x(t)\cos y(t) &\sinh x(t)\sin y(t)\\
0          &\cosh x(t) &- \sinh x(t)\sin y(t) &\sinh x(t)\cos y(t)\\
 \sinh x(t)\cos y(t) &-\sinh x(t)\sin y(t) &\cosh x(t) &0         \\
 \sinh x(t)\sin y(t) & \sinh x(t)\cos y(t) &0          &\cosh x(t)\\
\end{array}\right]
}
\end{align*}
with $x(t)\geq 0$ so that $p(\wt\gamma(t))=\gamma(t)$ ($\wt\gamma$
is a lift of $\gamma$), and let {\smaller
\[
\omega(t)= \left[
\begin{array}{rrrr}
\cos z(t)  &  -\sin z(t) & 0& 0  \\
\sin z(t)  & \cos z(t)   & 0 & 0\\
 0& 0 & \cos z(t) & \sin z(t) \\
 0& 0 & -\sin z(t) & \cos z(t)  \\
\end{array}\right].
\] }
 Put
$ \eta(t)=\wt\gamma(t)\cdot\omega(t). $ Then still
$p(\eta(t))=\gamma(t)$, and $\eta$ is another lift of $\gamma$. We
wish $\eta$ to be the horizontal lift of $\gamma$. That is, we
want $\eta'(t)$ to be orthogonal to the fiber at $\eta(t)$. The
condition is that
$\angles{\eta'(t),(\ell_{\eta(t)})_{*}(e_{3})}=0$, or
equivalently, $\angles{(\ell_{\eta(t)\inv})_*\eta'(t),e_3}=0$.
That is,
$$
\eta(t)^{-1} \cdot \eta'(t)= \alpha_{1}e_1  + \alpha_{2}e_{2}
$$
for some $\alpha_1,\alpha_2 \in\bbr$. From this, we get the
following equation:
\begin{equation}
\label{hyper-z} z'(t)=\sinh^2 x(t) y'(t).
\end{equation}

By virtue of Lemma \ref{rect-to-gen}, it will be enough to prove
the statement for a particular type of curves as follows:
Suppose we are given a rectangular region in the $xy$-plane
\begin{align*}
p \leq x \leq p+a, \quad q \leq y \leq q+b.
\end{align*}
Consider the image $R$ of this rectangle in  $\bbc H^1$ by the map
$$
(x,y)\mapsto \bfr(x,y)=(\cosh 2x, (\sinh 2x)(\cos y), (\sinh
2x)(\sin y)),
$$
$\bbc H^1$  with the $(+--)$ metric. The area of $R$ can be
calculated as follows:
$$
\bfr_x\x\bfr_y=((2\cosh 2x)(\sinh 2x), -(2\sinh^2 2x)(\cos y),
-(2\sinh^2 2x)(\sin y)).
$$
Now
\begin{align*}
||\bfr_x\x\bfr_y||
&=2|\sinh 2x|,\quad (+--)\text{-norm}\\
&=2 \sinh 2x,\quad\text{(because $x\geq 0$)}.
\end{align*}
Thus, the area is
\begin{align*}
\int_q^{q+b}\int_p^{p+a} 2\sinh 2x\ dx dy &=\left[\left[2 \sinh^2
x\right]_p^{p+a}\right]_q^{q+b} =2 b(\sinh^2(p+a)-\sinh^2(p)).
\end{align*}

On the other hand, the change of $z(t)$ along the boundary of this
region can be calculated using  condition (\ref{hyper-z}). Label
the  four vertices by $A(p,q)$, $B(p+a,q)$, $C(p+a,q+b)$, and
$D(p,q+b)$. $AB$ can be parametrized by $x(t)=p+at$, $y(t)=q$, $t
\in [0,1]$  so that $y'(t)=0$. For $BC$, $x(t)=p+a$, $y(t)=q+bt$,
$t \in [0,1]$. Then
$$
z(1)-z(0) =\int_0^1 z'(t) dt = \int\sinh^2 (p+a) b\ dt =
b\cdot\sinh^2 (p+a).
$$
Similarly, $z(t)$ does not change along $CD$, but on $DA$,
$x(t)=p$, $y(t)=q+b-bt$, $t \in [0, 1]$.  So
$$
z(1)-z(0) =\int_0^1 z'(t) dt = \int\sinh^2 (p) (-b) dt =-
b\cdot\sinh^2 (p).
$$
Thus the total vertical change of $z$-values, $z(1)-z(0)$, along
the perimeter of this rectangle is $ b\cdot(\sinh^2 (p+a)-\sinh^2
(p)) $ which is $1/2 $  times  the area.
\end{proof}
\bigskip

Now we turn to the general case
\[ S^1 \rightarrow \son \stackrel{p}{\rightarrow}  \bbc H^n. \]
We are viewing $ \son\cong U(1,n)/U(n)$, and $\bbc H^n \cong
U(1,n)/(U(1) \times U(n)). $ The Lie algebra of $U(1,n)$ is
$\fraku(1,n)$, and has the following canonical decomposition: $
\frakg=\frakh +\frakm, $ where
$$
\frakh=\fraku(1)+\fraku(n)= {\smaller  \left\{
\left[\begin{array}{cc} \lambda  & 0\\ 0 & B \end{array}\right]\
:\ \lambda + \bar{\lambda }=0,\  B\in \fraku(n)\right\} }
$$
and
$$
\frakm= {\smaller \left\{\left[\begin{array}{cc}
 0 & \bar{\xi}^{t} \\ \xi &0
\end{array} \right]\ :\ \xi  \in  {\bbc}^{n} \right\}.}
$$

\smallskip
\begin{lemma}
\label{geod-cond-hyper} A 2-dimensional subspace $\frakm'$ of
$\frakm\subset {\mathfrak u}(1,n)$ gives rise to a complete
totally geodesic submanifold of $\bbc H^n$ if and only if either
\begin{enumerate}
\item
$\frakm'$ is $J$-invariant (i.e., has a complex structure), or
\item
$\frakm'$ has tangent vectors $\bfv$ and $\bfw$ such that $\ol
\bfv \bfw- \bfv \ol \bfw=0$.
\end{enumerate}
Furthermore, for each of these cases, the pullback of the bundle
$S^1\ra \son\ra \bbc H^n$ by the inclusion is
 isomorphic to the standard bundle $S^1\ra SU(1,1)\ra \bbc H^1$ for {\rm (1)},
or the product bundle $S^1\x \bbc H^1$ for {\rm (2)},
respectively.
\end{lemma}

\begin{proof}
With the notation $m$ as above, let  $\bfv$ and $\bfw$  be
elements of $\frakm$ whose $\xi$'s are given by
 $$
\left[\begin{array}{c} x_1+iy_1\\ x_2+iy_2\\ \cdots\\
x_n+iy_n\end{array}\right] \quad\text{ for\, $\bfv$\, and}\quad
\left[\begin{array}{c} a_1+ib_1\\ a_2+ib_2\\ \cdots\\
a_n+ib_n\end{array}\right] \quad\text{ for\,  $\bfw$}.
$$
Then by Proposition  \ref{tot-geod-prop}, $\frakm'$ is a totally
geodesic sub-manifold if and only if
$[[\frakm',\frakm'],\frakm']\subset \frakm'$. Some calculations
show the following equality
$$
[[\bfv,\bfw],\bfv] = \sum_{k=1}^n (x_k a_k + y_k b_k) \bfv -
\sum_{k=1}^n (x_k^2 + y_k^2) \bfw - 3 \sum_{k=1}^n (x_k b_k - y_k
a_k) (i\bfv)
$$
holds. Therefore, $[[\bfv,\bfw],\bfv]=p \bfv + q\bfw$ for some
real $p$ and $q$ if and only if $i\bfv=p \bfv + q\bfw$ has
solution for some real $p$ and $q$.

Suppose $\frakm'$ has a complex structure. Then we can take $\bfv$
and $\bfw$   in  $\frakm'$ so that  $i\bfv=\bfw$ (so $a_k=-y_k$
and $b_k=x_k$ for all $k=1,2,\cdots,n$). Thus,
$[[\bfv,\bfw],\bfv]=p \bfv + q\bfw$ has a solution for $p$ and
$q$. Suppose $\ol \bfv \bfw-\bfv \ol \bfw=2\text{Im}(\ol \bfv
\bfw)=2\sum_{k=1}^n (x_k b_k - y_k a_k)=0$. Then clearly
$[[\bfv,\bfw],\bfv]=p \bfv + q\bfw$ has a solution for $p$ and
$q$.

Conversely, suppose $[[\bfv,\bfw],\bfv]=p \bfv + q\bfw$ has a
solution for $p$ and $q$. Then $i\bfv=p \bfv + q\bfw$ must have a
real solution for $p$ and $q$. Suppose $\sum_{k=1}^n (x_k b_k -
y_k a_k)\not=0$. Then, at least one of the summands is non-zero,
say $x_1 b_1 - y_1 a_1\not=0$. This means that we can find a new
basis for $\text{span}\{\bfv,\bfw\}$ with
$$
x_1=1,\ y_1=0;\quad a_1=0,\ b_1=1.
$$
Then the equation $ p \bfv + q \bfw = i\bfv $ is quickly reduced
to $p=0$ and $q=1$ (from $k=1$), and hence we obtain $
x_k=b_k,\quad y_k=-a_k $ for all $k=2,\cdots,n$. This shows
$\bfw=i\bfv$, and the space spanned by $\bfv$ and $\bfw$ has a
complex structure.

For the second part of the statement, it is enough to observe that

\centerline {$ [\bfv,\bfw]= \left[\begin{array}{cc}
 \lambda & 0 \\ 0 &0
\end{array} \right]
$}

\noindent  where $\lambda=\ol \bfv \bfw - \bfv\ol \bfw$. If
$\lambda=0$, then the distribution $\frakm'$ is integrable, and
the bundle is trivial.
\end{proof}

By combining Theorem \ref{area-u11} and Lemma
\ref{geod-cond-hyper}, we have now

\begin{theorem}
\label{thm-hyper} Let $S^1\ra  \son \ra \bbc H^n$ be the natural
fibration. Let $S$ be a complete totally geodesic 2-dimensional
surface in $\bbc H^n$, and $\xi_S$ be the pullback bundle over
$S$. Let $\gamma$ be a piecewise smooth, simple closed curve on
$S$. Then the holonomy displacement along $\gamma$ is given by
$$
V(\gamma)=e^{\frac{1}{2} A(\gamma) i}\text{\ or\ } e^{0 i}\ \in
S^1
$$
where $A(\gamma)$ is the area of the region on the surface $S$
surrounded by $\gamma$, depending on whether $S$ is a complex
submanifold or not.
\end{theorem}

\smallskip
Since the length of $\eta$ is half of the length of $\gamma$, we
have
\begin{corollary}
Suppose $\gamma$ is a piecewise smooth, simple closed curve
parametrized by arc length. Then the Hopf torus in $\son $ over
$\gamma$ is isometric to the torus generated by the lattice
$\{(2\pi,0),(\, A(\gamma)/2 \, , \, L(\gamma)/2 \, )\}$  in
$\bbr^2$, where $L(\gamma)$ is the length of $\gamma$.
\end{corollary}
\smallskip

\section{The complex Heisenberg group $\NILn$}

We consider $\NILn$. This is  $\bbr\x\bbc^n$ with group operation
given by
$$
(s,\bz)(t,\bz')=(s+t+2\ \im\{\ol\bz\bz'\},\ \bz+\bz'),
$$
where $\im\{\ol\bz\bz'\}$ is the imaginary part of the complex
number $\ol z_1 z'_1+ \ol z_2 z'_2+\cdots +\ol z_{n} z'_{n}$ for
$\bz=(z_1,z_2,\cdots,z_{n})$, $\bz'=(z'_1,z'_2,\cdots,z'_{n}) \in
{\bbc}^{n}$. This is a 2-step nilpotent Lie group with center
$\calz(\NILn )=\bbr$. In the case of $n=1$, $\NIL1$ is isomorphic
to the ordinary 3-dimensional Heisenberg group by
$$
(4z-2xy,x+iy)\longleftrightarrow \left[\begin{array}{ccc}
1 & x & z \\
0 & 1 & y   \\
0 & 0 & 1 \\
\end{array}\right].
$$

For the sake of computations, we use the following affine
representation of $\NILn$ into $\aff(2n+1)\subset\GL(2n+2)$:
\small{
\[
\left(s, \left[\begin{array}{c}
x_1+i y_1\\
x_2+i y_2\\
\cdots\\
\cdots\\
x_n+i y_n\\
\end{array}\right]\right)
\lra \left(
\begin{array}{cccccccc}
1 & -2y_{1} & 2x_{1}& \cdots  &-2y_{n} & 2x_{n}& s \\
0 & 1 & 0&   \cdots &  & 0 & 0 &  x_{1} \\
0 & 0 & 1&   \cdots  & 0 & 0 &  y_{1}\\
\cdot & \cdot & \cdot  & \cdots & \cdot  & \cdot &
\cdot \\
\cdot & \cdot & \cdot  & \cdots & \cdot  & \cdot &
\cdot \\
0 & 0 & 0&  \cdots   & 1 & 0 &  x_{n}\\
0 & 0 & 0&   \cdots   & 0 & 1 &  y_{n}\\
0 & 0 & 0&   \cdots   & 0 & 0 &  1\\
\end{array}\right)
\]
The Lie algebra  has a  following orthonormal basis {\tiny
\begin{eqnarray}
e_{1}= \left(
\begin{array}{cccccccc}
0 & 0 & 2& \cdots  &0 & 0& 0 \\
0 & 0 & 0&   \cdots   & 0 & 0 &  1 \\
0 & 0 & 0&   \cdots  & 0 & 0 &  0\\
\cdot & \cdot & \cdot  & \cdots & \cdot  & \cdot &
\cdot \\
\cdot & \cdot & \cdot & \cdots & \cdot  & \cdot &
\cdot \\
0 & 0 & 0&  \cdots  & 0 & 0 &  0\\
0 & 0 & 0&   \cdots   & 0 & 0 &  0\\
0 & 0 & 0&   \cdots   & 0 & 0 &  0\\
\end{array}\right)
,  \nonumber \;\; & e_{2}= \left(
\begin{array}{ccccccccc}
0 & -2 & 0& \cdots  &0 & 0& 0 \\
0 & 0 & 0&   \cdots   & 0 & 0 &  0 \\
0 & 0 & 0&   \cdots  & 0 & 0 &  1\\
\cdot & \cdot & \cdot  & \cdots & \cdot  & \cdot &
\cdot \\
\cdot & \cdot & \cdot  & \cdots & \cdot  & \cdot &
\cdot \\
0 & 0 & 0&  \cdots  & 0 & 0 &  0\\
0 & 0 & 0&   \cdots   & 0 & 0 &  0\\
0 & 0 & 0&   \cdots  & 0 & 0 &  0\\
\end{array}\right),\;\; \\
\vdots  \;\;\;\qquad \qquad  \qquad \qquad & \qquad  \vdots \nonumber  \\
e_{2n}=\left(
\begin{array}{ccccccccc}
0 & 0 & 0& \cdots  &-2 & 0& 0\\
0 & 0 & 0&   \cdots   & 0 & 0 &  0 \\
0 & 0 & 0&   \cdots   & 0 & 0 &  0\\
\cdot & \cdot & \cdot  & \cdots & \cdot  & \cdot &
\cdot \\
\cdot & \cdot & \cdot  & \cdots & \cdot  & \cdot &
\cdot \\
0 & 0 & 0&  \cdots   & 0 & 0 &  0\\
0 & 0 & 0&   \cdots   & 0 & 0 &  1\\
0 & 0 & 0&   \cdots & 0 & 0 &  0\\
\end{array}\right) &
e_{2n+1}=\left(
\begin{array}{ccccccccc}
0 & 0 & 0& \cdots &0 & 0& 1\\
0 & 0 & 0&   \cdots   & 0 & 0 &  0 \\
0 & 0 & 0&   \cdots   & 0 & 0 &  0\\
\cdot & \cdot & \cdot & \cdots & \cdot  & \cdot &
\cdot \\
\cdot & \cdot & \cdot & \cdots & \cdot  & \cdot &
\cdot \\
0 & 0 & 0&  \cdots  & 0 & 0 &  0\\
0 & 0 & 0&   \cdots   & 0 & 0 &  0\\
0 & 0 & 0&   \cdots   & 0 & 0 &  0\\
\end{array}\right)  \nonumber
\end{eqnarray} }

\noindent which defines a left-invariant Riemannian metric on
$\NILn$. The short exact sequence of groups
$$
1\ra\bbr\ra\NILn\stackrel{p}\ra\bbc^n\ra 1
$$
is a fiber bundle, which is topologically trivial. The left
invariant metric naturally induces a connection on this principal
$\bbr$-bundle. There is a unique metric on $\bbc^n$ (standard
Euclidean metric) which makes the projection map $p$ a Riemannian
submersion.

\bigskip

\begin{theorem}
\label{thm-nil} Let $1\ra\bbr\ra\NILn\stackrel{p}\ra\bbc^n\ra 1$
be the central short exact sequence of the complex Heisenberg
group. Let $S$ be a complete totally geodesic plane in $\bbc^n$,
and $\xi_S$ be the pullback bundle over $S$. Let $\gamma$ be a
piecewise smooth, simple closed curve on $S$. Then
$$
V(\gamma)= e(\xi_S)\cdot A(\gamma)
$$
where $A(\gamma)$ is the area of the region on the surface $S$
surrounded by $\gamma$, and the number $e(\xi_S)$ is determined by
the equality $[\bfv,\bfw]=e(\xi_S) e_{2n+1}$ for an orthonormal
basis $\{\bfv,\bfw\}$ for the tangent space of $S$.
\end{theorem}

\begin{proof}
Every complete totally geodesic submanifold of $\bbc^n$ is an
$\bbr$-linear subspace of $\bbc^n$. Therefore
$S=\text{span}\{\bfv,\bfw\}$ for some  orthonormal basis
$\bfv,\bfw $ where $\bfv=\sum_{j=1}^{n}(a_j+ib_j)$ and
$\bfw=\sum_{j=1}^{n}(c_j+id_j)\in\bbc^n$.
Then $\gamma$ is of the form
$$
\gamma(t)=x(t)\bfv + y(t)\bfw\in S\subset\bbc^n,
$$
where $x(t)$ and $y(t)$ are scalars. We want to find a curve
$z(t)$ in $\bbr$ so that $ \eta(t)=(z(t),\gamma(t)) $ is
orthogonal to the fiber for every $t$. In other words,
$$
\angles{\eta'(t),(\ell_{\eta (t)})_{*} e_{2n+1}}=0,
$$
where $\ell$ is the left translation. This is equivalent to
$\angles{(\ell_{\eta(t)}\inv)_*\eta'(t), e_{2n+1}}=0.$ Note that
$\eta(t)\inv = (-z(t)+2\
\im\{\overline{\gamma(t)}\gamma(t)\},-\gamma(t))$. Using the
affine representation, $(\ell_{\eta (t)^{-1} })_{*}\eta'(t)$ is
$$
{\smaller \left(
\begin{array}{ccccc}
0 & -2( x'(t)b_{1} + y'(t) d_{1} )
 &  \cdots &  z'(t) - 2 (  x(t) y'(t) -   x'(t) y(t))\ \im\{\ol\bfv  \bfw\} \\
0 & 0 &  \cdots  &  x'(t)a_{1} + y'(t) c_{1}   \\
0 & 0 &    \cdots  &  x'(t)b_{1} + y'(t) d_{1}  \\
\cdot &  \cdot & \cdots  &
\cdot \\
\cdot & \cdot & \cdots  &
\cdot \\
0 & 0 &   \cdots  & x'(t)a_{n} + y'(t) c_{n}  \\
0 & 0 &    \cdots  & x'(t)b_{n} + y'(t) d_{n} \\
0 & 0 &    \cdots  &  0\\
\end{array}\right)}
$$
where $\im\{  \ol\bfv \bfw \}=\sum_{j=1}^{n}(a_{j}d_{j} -
c_{j}b_{j})$. Note that $\im \{\ol\bfv \bfw \}=\im\{\ol\bfv'
\bfw'\}$ for any  orthonormal basis $\{ \bfv',\bfw' \}$.
  The equation $\angles{(\ell_{\eta (t)^{-1}
})_{*}\eta'(t),e_{2n+1}}=0$ gives rise to
\begin{equation}
\label{nil-z} z'(t) - 2 (x(t)y'(t) - x'(t)y(t))\ \im\{\ol\bfv \bfw
\}=0.
\end{equation}
Suppose we are given a rectangular region on $xy$-plane
\begin{align*}
p \leq x \leq p+a,\quad q \leq y \leq q+b.
\end{align*}
Consider the image $R$ of this rectangle in $S\subset\bbc^n$ by
the map
$$
(x,y)\mapsto x\bfv+y\bfw.
$$

\noindent Then $R$ is a rectangle with vertices $ p\bfv+q\bfw$,
$(p+a)\bfv+q\bfw$, $(p+a)\bfv+(q+b)\bfw$, $p\bfv+(q+b)\bfw$. Let
$\gamma(t)$ be the piecewise linear boundary curve. It can be
represented by $((p+4at)\bfv,q\bfw) $        for   $0\leq t \leq
1/4$, $((p+a)\bfv,(q+b(4t-1))\bfw)$    for   $1/4\leq t \leq 1/2$,
$((p+a(3-4t))\bfv,(q+b)\bfw) $   for $ 1/2\leq t \leq 3/4$,
$(p\bfv,(q+b(4-4t))\bfw)$    for   $3/4\leq t \leq 1 $.

\noindent Then, from equation (\ref{nil-z}),
\[
z( 1 )-z( 0 ) = 2\int_{0}^{1} (x(t) y'(t) -   x'(t) y(t))\ \im\{
\ol\bfv \bfw \}  dt = 4 ab\ \im\{ \ol\bfv \bfw \}.
\]
On the other hand,
\[ [\bfv,\bfw] = 4 \left(
\begin{array}{cccccccc}
0 & 0 & 0& \cdots  &0 & 0&  \im\{\ol\bfv \bfw \} \\
0 & 0 & 0&   \cdots   & 0 & 0 &  0 \\
0 & 0 & 0&   \cdots   & 0 & 0 &  0\\
\cdot & \cdot & \cdot  & \cdots & \cdot  & \cdot &
\cdot \\
\cdot & \cdot & \cdot  & \cdots & \cdot  & \cdot &
\cdot \\
0 & 0 & 0&  \cdots   & 0 & 0 &  0\\
0 & 0 & 0&   \cdots   & 0 & 0 &  0\\
0 & 0 & 0&   \cdots   & 0 & 0 &  0\\
\end{array}\right). \\ \]
This means that $[\bfv,\bfw] = 4\, \im\{ \ol\bfv \bfw \}
e_{2n+1}=e(\xi_S)\ e_{2n+1}$ so that $V(\gamma)= e(\xi_S)\cdot
A(\gamma)$ with $e(\xi_S)=4\, \im\{ \ol\bfv \bfw \}$.

Having shown the statement for rectangular regions, now we apply
Lemma \ref{rect-to-gen} to conclude that the same formula holds
for any piecewise smooth, simple closed curve.
\end{proof}

\begin{corollary}
Suppose $\gamma$ is  a piecewise smooth, simple closed curve
parametrized by arc-length. Then the Hopf cylinder in $\NILn$ over
$\gamma$ is isometric to the cylinder generated by the translation
$(x,y) \mapsto (x+e(\xi_S)\, A(\gamma) , \, y+L(\gamma)\, )$ on
$\bbr^2$, where $L(\gamma)$ is the length of $\gamma$.
\end{corollary}

\bigskip

\noindent \textbf{ Acknowledgements.}
\smallskip
The first  author was supported by  KOSEF R01-2004-000-10183-0.
This work was done while the first author was visiting  University
of Oklahoma. He wishes to express his sincere thanks for the
hospitality.

\medskip

\bibliographystyle{amsplain}

\end{document}